\documentstyle[]{article}
  \begin{document}
  \title{Cartan Pairs \footnote{To the memory of Professor Jan Rzewuski.}}
\author{
Andrzej Borowiec \thanks{Supported by the State Research Committee 
KBN No 2 P302 023 07.}\\Institute of Theoretical Physics,
University of Wroc{\l}aw\\ 
Pl. Maxa Borna 9, 50-204 Wroc{\l}aw, Poland\\ 
e-mail: borow@ift.uni.wroc.pl}
\maketitle
\def\ba{\begin{array}}
\def\ea{\end{array}}
\def\lra{\longrightarrow}
\def\ra{\rightarrow}
\def\ld{\ldots}
\def\cd{\cdot}
\def\be{\begin{equation}}
\def\ee{\end{equation}}
\def\bm{\begin{em}}
\def\en{\end{em}}
\def\bk{I \kern-.25em k} 
\def\lin{\hbox{lin}\,}
\def\der{\hbox{der}\,}
\def\derk{\der_{\bk}\,}
\def\bimod{\hbox{bimod}\,}
\def\im{\hbox{im}\,}
\def\ker{\hbox{ker}\,}
\def\ann{\hbox{ann}\,}
\def\bot{\,_b\!\otimes}
\def\nd{\hbox{End}\,}
\def\hom{\hbox{Hom}\,}
\def\alg{\hbox{alg}\,}

\def\id{\hbox{id}}
\def\idA{\id_A}
\def\idV{\id_V}
\def\ot{\otimes}
\def\c{\circ}
\def\bot{\,_b\!\otimes}
\def\l{\lambda}
\def\lc#1{#1^{\lambda}}
\def\p{\partial}
\def\rp#1{#1^{\p}}
\def\lp#1{^{\p}\!#1}
\def\M{\Omega}
\begin{abstract}
A new notion of Cartan pairs as a substitute of notion of vector fields
in the noncommutative geometry is proposed.
The correspondence between Cartan pairs and differential calculi 
is established.
\end{abstract}
\section{Introduction.}

A big part of the
classical differential geometry on manifold $\M$, see {\it e.g.} \cite{W}, 
belongs to the theory of modules over a commutative algebra
${\cal F}\M$ of smooth scalar valued functions on $\M$.
One defines a tangent $T\M$ and cotangent $T^{*}\M$ vector bundles. Their 
sections, vector fields and one forms respectively, constitute modules 
${\cal X}\M$ and $\Lambda^ 1 \M$ over ${\cal F}\M$ and are basic 
differential--geometric objects on $\M$. Both notions (of vector fields and
that of one forms) enter the game on equal rights and are mutually dual.
In particular, $\Lambda^1 \M$ can be identify with a module of ${\cal 
F}\M$--linear mappings from ${\cal X}\M$ into ${\cal F}\M$ and evaluation of
one form $\omega$ on a vector field $X$ provides a ${\cal F}\M$--bilinear
pairing $\omega(X)\equiv <X, \omega>\in {\cal F}\M$ between these modules.
Also an action of vector fields on functions and an external differentiation
of functions are dual each other via famous Cartan formulae
$$
X(f)\equiv <X,\, df>\equiv \hbox{i}_X df \in {\cal F}\M \ . \eqno (1.1)
$$
It appears that the Leibniz rule
$$
d(fg) = (df).\,g + f.\,dg  \eqno (1.2)
$$
for an external differential $d: {\cal F}\M\ra \Lambda^1 \M $ of functions
into one forms
and derivation property of vector fields (the Leibniz rule for an
"internal" derivation $X:{\cal F}\M\ra {\cal F}\M $)
$$
X(f\,g) = X(f)\,g + f\,X(g) \eqno (1.3)
$$
are related each other by (1.1). A vector field can be 
alternatively defined as a  derivation of ${\cal F}\M$ {\it i.e.}
as an endomorphism of ${\cal F}\M$ satisfying (1.3). Therefore, the
module of vector fields bears a Lie algebra structure. 

To be precise, one should distinguish
between a vector field $X$ as a smooth section of $T\M$ and its isomorphic 
image $X \in \hbox{Der}({\cal F}\M)\subset \hbox{End}({\cal F}\M)$, 
which acts on functions via (1.1). These ideas have served as a
basis for an algebraic generalization of concept of vector fields,
so called  {\it Lie--Cartan pairs} \cite{KS} or {\it Lie pseudoalgebras} 
\cite{Mc}, see also \cite{JK} for supersymmetric generalization and 
\cite{Mc} for overview and historical remarks. An attempt to generalize 
this concept to noncommutative case within the framework of braided Lie 
algebras \cite{G}  was performed in \cite{Wl}
(c.f. also \cite{P} and \cite{OPR}). The notion of Lie algebras of vector 
fields for quantum groups has been introduced by Woronowicz \cite{Wor}.

Passing to the noncommutative case, the duality between 
forms and vector fields fails. 
Unlike in the commutative case, the noncommutative differential 
calculus is developed  mainly in the covariant approach ({\it i.e.} by
means of differential forms). A satisfactory concept
of noncommutative vector fields has not been formulated yet.
The reason is that the Leibniz rule (1.2) for an external differential
remains unchanged also in the noncommutative setting while that one for 
vector fields (1.3) has to be modified.
The aim of the present note is to fill this gap. We propose a new notion of
{\it Cartan pairs} as a substitute for a concept of vector fields. Our
approach is similar to the Lie--Cartan approach but we have no analogue
of Lie bracket. We explore a bimodule structure instead. A Cartan pair
consist an $\bk$--algebra $A$ ($\bk$ being a commutative ring) and 
$A$--bimodule $M$ with suitable action of $M$ on $A$. We show that a dual
object to a Cartan pair is a differential calculi on an
algebra $A$. Our main result is that (1.1) allows
to reconstruct the "external" differential if we are given an action
of generalized vector fields  and conversely to find out the 
action by means of differential. An example of such action for a given
noncommutative calculus can be found in \cite{BK} (c.f. \cite{DMH}).

Henceforth $\bk$ denotes some fixed unital and commutative ring. 
Algebras are unital associative $\bk$-algebras and homomorpisms
are assumed to be unital. 
All objects considered here  are first of all $\bk$-modules. 
All maps are assumed to be $\bk$-linear maps. 

Let $M$ be an $(A, A)$--bimodule ($A$--bimodule in short). We shall
denote by dot "." the both: left and right multiplication by elements
from $A$. For example, by bimodule axioms, one has $(f.x).g = f.(x.g) = 
f.x.g$ for $f, g\in A$ and $x\in M$.

The present note
has a preliminary character. 
The full version of it with
more details and proofs will be published elsewhere.

\section{Cartan pairs.}
Let $A$ be an $\bk$--algebra and $M$ an $A$--bimodule. By an {\it action}
of $M$ on $A$ we mean a $\bk$--linear mapping 
$\beta\in\hom_{\bk}(M,\, \nd_{\bk}(A))$. We shall also write
$M\ni x\mapsto x^\beta\in\nd_{\bk}(A)$ or
$A\ni f\mapsto x^{\beta}(f)\in A$ to denote the action.

{\bf DEFINITION 2.1}. Let $\bk$ be a commutative and unitary ring, and let $A$
be an unitary $\bk$--algebra. A {\it right Cartan pair} over $\bk$
and $A$ is an $A$--bimodule $R$ together with a {\it right action}
$\rho :R\ra \nd_{\bk}(A)$, such that
$$
(f.X)^\rho (g)\ =\ f\, X^\rho (g) \eqno (2.1)$$
and
$$
X^\rho (f\, g)\ =\ X^\rho (f)\,g + (X.f)^\rho (g) \eqno (2.2)$$

Observe that in the case of commutative algebras $X.f=f.X$ (c.f. Remark 3.7
below) the formulae (2.1) and (2.2) set a generalization of the Leibniz 
rule (1.3) we have been looking for.

In a similar manner we define a {\it left Cartan pair} $(L, \l)$ consisting
a bimodule $L$ and its {\it left action} $L\ni X\mapsto X^\l \in\nd_{\bk}(A)$.
Now the properties (2.1) and (2.2) must be replaced by
$$
\lc{(X.g)}(f)\ =\ \lc{X}(f) \,g \eqno (2.3)$$
$$
\lc{X}(f\,g)\ =\ f\,\lc{X}(g) + \lc{(g.X)(f)}  \eqno (2.4)$$

{\bf DEFINITION 2.2}. A left (resp. right) Cartan pair $(M, A)$ is called a 
bimodule of left (resp. right) {\it generalized vector fields} on $A$ if the 
corresponding action $\lambda$ (resp. $\rho$) is faithful.

{\bf EXAMPLE 2.3}. Let $\M$ be a manifold and $\bk$ be a field of real 
numbers. Take $A={\cal F}\M$ and $M={\cal X}\M$ together with a canonical 
action of vector fields on function via derivations. Since algebra is 
commutative, the module ${\cal X}\M$ can be considered as a bimodule with a 
left and right multiplication coinciding. Then $({\cal X}\M, {\cal F}\M)$
is at the same time a left and right Cartan pair. Of course it is 
a bimodule of generalized vector fields.

\section{Dual of bimodule.}

Let  $R$ be a right $A$--module. Recall (see {\it e.g.} \cite{Bou}) 
that $R^*$ dual of $R$ is defined
as a collection of all right $A$--module maps from $R$ int $A$, {\it i.e.}
$R^* =\hom_A(R, A)$. For every ordered pair of elements $x\in R$ and
$X \in R^*$, the element $X(x)\in A$, the evaluation of $X$ on $x$
is denoted by $<X,\, x>$.  $R^*$ bears
a canonical left $A$--module structure, therefore $<,>:R^*\times R\ra A$
defines the canonical $A$--bilinear form ({\it pairing}).
Summing up the following relations hold true
$$
<X,\, x+y> \ =\  <X,\, x>+<X,\, y>  \eqno (3.1)$$
$$
<X,\, x.f> \ =\  <X,\, x>\, f \eqno (3.2)$$
$$
<X+Y,\, x> \ =\  <X,\, x>+<Y,\, x> \eqno (3.3)$$
$$
<f.\,X,\, x> \ =\  f\, <X,\, x> \eqno (3.4)$$
where, $X,\,Y\in R$ and $f\in A$.

For a right $A$--modules map $\alpha\in \hom_A(R_1,\,R_2)$ one defines
its transpose $\alpha^T$ as a left $A$--module
map $\alpha^T\in\hom_A(R^*_2,\, R^*_1)$ by the formulae \cite{Bou}
$$
<\alpha^T(X_2),\, x_1>_1 \ =\  <X_2,\, \alpha(x_1)>_2$$
where, $x_i\in R_i$ and $X_i\in R^*_i$, $i = 1,2$.

Let now $M$ be an $A$--bimodule and let $M^* \ =\  \hom_{(-,A)}(M,\, A)$ 
denotes its right dual, {\it i.e.} dual of $M$ as a right $A$--module. 
For any element $f\in A$ left multiplication by $f$ is a right module map 
$f. \in\hom_{(-,A)}(M,\, M)$. It is easy to check that its transpose
$(f.\,)^T\equiv \,.f$ is a right multiplication in $M^*$ and that 
with this multiplication $M^*$ becomes a bimodule.

{\bf DEFINITION 3.1}. The $A$--bimodule $M^* = \hom_{(-,A)}(M,\, A)$ with 
the canonical left module structure (3.3), (3.4) and with transpose right 
multiplication
$$
<X.f,\, x> \ =\  <X,\, f.x> \eqno (3.5)$$
is called a right dual of a bimodule $M$.

In a similar way one defines a {\it left dual} $^*\!M = \hom_{(A,-)}(M,\,A)$ 
of bimodule $M$ with a canonical left and transpose right $A$--module 
structure. In this case one has
$$
<x+y,\, X> \ =\  <x,\, X>+<y,\, X> \eqno (3.6)$$
$$
<f.\,x,\, X> \ =\  f\, <x,\, X> \eqno (3.7)$$
$$
<x,\, X+Y> \ =\  <x,\, X>+<x,\, Y> \eqno (3.8)$$
$$
<x,\, X.\,f> \ =\   <x,\, X>\,f \eqno (3.9)$$
$$
<x.\,f,\, X> \ =\  <x,\, f.\,X> \ \ . \eqno (3.10)$$

It is interesting to compare a left dual of a right dual of a bimodule
$M$ with $M$.

{\bf PROPOSITION 3.2}. {\em There is a canonical $A$--bimodule map from 
$M$ into $^*\!(M^*)$ 
(resp. $(^*\!M)^*$) $x\mapsto \tilde x$  given by the formulae
$$
<X,\, \tilde x> \ =\  <X,\, x>\ \ \ 
(\hbox{resp.}\ <\tilde x,\, X> \ =\  <x,\, X>\ )\ . \eqno (3.11)$$
In general, it is neither injective nor surjective .}

{\bf DEFINITION 3.3}. A bimodule $M$ is called right (resp. left) reflexive if
$^*\!(M^*)\equiv M$ (resp. $(^*\!M)^*\equiv M$) {\it i.e.} when the
corresponding canonical map (3.11) is a bimodule isomorphism.

{\bf DEFINITION 3.4}. An $A$--bimodule $M$ is called a right (resp. left) free
$A$--bimodule if it is so as a right (resp. left) $A$--module.

{\bf DEFINITION 3.5}. An $A$--bimodule $M$ is called a right (resp. left) 
finitely generated
bimodule if it is so as a right (resp. left) $A$--module.

{\bf LEMMA 3.6}. {\em Let $M$ be a right (resp. left) free and finitely 
generated $A$--bimodule. Then a right (resp. left) dual of $M$
is a left (resp. right) free finitely generated bimodule.
Moreover, $M$ is a right (resp. left) reflexive.}

{\bf REMARK 3.7}. Assume that $A$ is commutative. Each $A$--module becomes 
automatically an $A$--bimodule with the same left and right multiplication. 
In this case the
three notions of dual, namely: left and right dual of bimodule and
dual of module over commutative algebra, coincide. Of course, the module
of the classical vector fields ${\cal F}\M$ over manifold $\M$ is 
reflexive (c.f. Example 2.3).

\section{Main results.}

It appears that differential calculi investigated recently by many
authors in the context of quantum groups and noncommutative 
geometry are nothing but derivations of an algebra with values in
a bimodule (c.f. \cite{BKO,BK,DMH,PW,WZ}).
Recall that a $\bk$-derivation $d$ of $A$ to $M$,
$d\in\hbox{Der}_{\bk}(A, M)$, is a $\bk$-linear mapping from
$A$ into $M$ such that the Leibniz rule (1.2) is satisfied.
The pair $(M, d)$ is said to be {\it first order calculus} or 
{\it first order differential} on an algebra $A$ with values in 
an $A$--bimodule $M$ or shortly {\it $M$-valued calculus} on $A$.
Each $\bk$-derivation vanishes on scalars from $\bk$.

Let now $(M, d)$ be a calculus on an algebra $A$. The differential
$d$ and formulae (1.1) defines an action of the right dual $M^*$ on $A$. 
This action 
$$A\ni f\mapsto\  \rp X(f)\ \equiv\  <X,\, d f> \eqno (4.1)$$
will be called a {\it right partial derivatives} along the element $X\in 
M^*$ with respect to the calculus $(M, d)$. One  uses $\rp X$ instead of 
more traditional notation $\p / \p X$. It can be checked that this action
satisfies axioms of right Cartan pair. Therefore, to each differential
calculus $(M, d)$ on $A$ we can  associate a unique right 
Cartan pair of right partial derivatives $(M^*, \p)$ of $(M, d)$. The 
converse statement is also true: to each right Cartan pair $(R, \rho)$ one 
can associate a unique differential calculus $(^*\!R, d_\rho)$ where,
$d_\rho:A\ra\ ^*\!R$  is defined by formulae (4.2) below. Thus we have

{\bf MAIN THEOREM}. {\em
Let $(M, d)$ be a calculus on $A$. Then $M^*$ together with an action (4.1),
via the right partial derivatives, becomes a a right Cartan pair 
$(M^*, \p)$ on $A$. Moreover, if the module $M$ of one forms is spanned by 
differential ({\it i.e.} $M=A.dA$) then the action $\p$ is faithful.

Conversely, let $(R, \rho)$ be a right Cartan pair on $A$. Then the formulae
$$
<X, d_\rho f>\  =\  X^\rho (f) \eqno (4.2)$$
for each $X\in R$, determines $d_\rho f$ as an element of a left dual $^*\!R$ 
of the bimodule $R$. The mapping $d_\rho : A\ra\ ^*\!R$ defines an 
$^*\!R$--valued calculus $(^*\!R, d_\rho)$ on $A$.

In a case of a right reflexive bimodule $M=^*\!(M^*)$ one has 
$d=d_\p$ and $\rho=\p_\rho$}.

In a similar way an action 
$$A\ni f\mapsto\  \lp X(f)\ \equiv\  <d f,\, X> \eqno (4.3)$$
determines a left Cartan pair structure on $^*\!M$ (left partial derivatives).

Therefore to each differential calculus one can canonically associate a 
right (resp. left) Cartan pair of partial derivatives. Conversely, for each
left (resp. right) Cartan pair there exists an associated differential 
calculus on an algebra $A$. In a case of reflexive bimodule a successive 
application of above canonical constructions give rise to the initial object.

An application of Cartan pairs in the theory of noncommutative vector
bundles and connections (c.f. \cite{DMH}) will be investigated elsewhere.

  \end{document}